\documentclass{article}

\usepackage[authoryear]{natbib}
\usepackage{url,doi,amssymb,amsmath,bussproofs,tikz,palatino,mathpazo}
\hypersetup{colorlinks,hidelinks}
\usetikzlibrary{automata,arrows}

\author{Richard Zach}
\title{Logic in Mathematics and Computer Science\thanks{To appear in Elke Brendel, Massimiliano Carrara, Filippo Ferrari, Ole Hjortland, Gil Sagi, Gila Sher, and Florian Steinberger, eds. \emph{The Oxford Handbook of Philosophy of Logic}. Oxford: Oxford University Press, 2024.}}

\begin{document}
\maketitle

\begin{abstract}
  Logic has pride of place in mathematics and its 20th century
  offshoot, computer science.  Modern symbolic logic was developed, in
  part, as a way to provide a formal framework for mathematics: Frege,
  Peano, Whitehead and Russell, as well as Hilbert developed systems
  of logic to formalize mathematics. These systems were meant to serve
  either as themselves foundational, or at least as formal analogs of
  mathematical reasoning amenable to mathematical study, e.g., in
  Hilbert's consistency program. Similar efforts continue, but have
  been expanded by the development of sophisticated methods to study
  the properties of such systems using proof and model theory. In
  parallel with this evolution of logical formalisms as tools for
  articulating mathematical theories (broadly speaking), much progress
  has been made in the quest for a mechanization of logical inference
  and the investigation of its theoretical limits, culminating
  recently in the development of new foundational frameworks for
  mathematics with sophisticated computer-assisted proof systems.  In
  addition, logical formalisms developed by logicians in mathematical
  and philosophical contexts have proved immensely useful in
  describing theories and systems of interest to computer scientists,
  and to some degree, vice versa. Three examples of the influence of
  logic in computer science are automated reasoning, computer
  verification, and type systems for programming languages.
\end{abstract}

\section{Introduction}

Modern logic got its start in two research programs, both intimately
tied to mathematics. The first was the mathematization of logic in the
work of Boole and the algebraic logicians of the 19th century. Boole
noticed that logical operations and relations, such as union,
intersection, and containment of concepts, and disjunction,
conjunction, and entailment of propositions, obey laws that can be
formulated as algebraic equations. The second was the formalization of
mathematical statements and of logical inference in the work of Frege,
Peano, Peirce, Whitehead, Russell, and Hilbert. Although their aims
and philosophical outlook diverged widely, they shared one fundamental
conviction. They all thought that in order to clarify the content of
mathematical statements, and to clarify fundamental concepts such as
mathematical inference, proof, and even (especially in the case of
Hilbert) mathematical existence, consistency, and independence of axioms,
it is necessary to \emph{formalize} mathematical theories. It was this
work that gave rise to the formalized systems of logic which
logicians now develop, expand, modify, and study, and with which the
philosophy of logic is concerned.

Languages for the formalization of mathematics were at first developed
purely syntactically, i.e., without a clear idea of how the symbols in
them were to be interpreted. But massive advances were made between
Frege's \emph{Begriffsschrift} \citep{Frege1879} and the first textbook
presentations of mathematical logic in the 1930s. Frege and Peirce
introduced polyadic predicates, propositional connectives, first- and
higher-order quantifiers, and identity. Whitehead and Russell
developed type theory, and Hilbert and others identified the
first-order fragment of classical logic. 

Frege, Whitehead, Russell, and Hilbert also provided axiomatizations
of their logical systems, leading to a clearly defined notion of
proof. Once a formal language and proof system were available, it
became natural to ask questions about this formal framework, and to
answer such questions with mathematical precision. This set the stage
for the development of model theory, soundness and completeness
theorems, decidability and undecidability, and the investigation of
specific mathematical theories.

The investigation of specific mathematical theories as formal
axiomatic systems is no doubt the most important and fundamental
contribution logic has made to mathematics. Some major early results
include the undecidability and incompleteness of axiomatized theories
of arithmetic, the formal axiomatization and investigation of set
theories, the consistency of the axiom of choice, and the decidability
of the theory of the real numbers. These results were made possible by
the development of formal logic, even if the results themselves are
mathematical and not, strictly speaking, logical results. That is,
they concern mathematical theories, and their proofs use mathematical
methods.

Until the 1920s, the logical systems introduced and investigated were
mostly classical: conditionals were material, excluded middle and double
negation elimination not questioned, truth values restricted to two,
and modalities and other intensional notions not considered. The only
questions that arose which one might now consider philosophical had to
do with higher types and quantification, e.g., whether impredicative
types should be allowed. In the 1920s, philosophers started to become
interested in the new symbolic logic, and philosophical questions gave
rise to non-classical variations: C. I. Lewis, motivated by the
paradoxes of the material conditional, began the study of modal logics
\citep{Lewis1918}. Łukasiewicz, motivated by the problem of future
contingents, introduced the first many-valued logics
\citep{Lukasiewicz1920}. It wasn't until the 1990s that the
philosophical development of non-classical logics came back to
mathematics. With few exceptions, mainstream mathematics has so far,
however, not been particularly interested in the development of
mathematical theories on the basis of many-valued, relevant, modal, or
paraconsistent logics.

A notable exception is intuitionistic logic. This is not surprising,
as intuitionistic logic arose out of a mathematical debate, namely the
foundational crisis of the 1920s. L. E. J. Brouwer proposed a
wholesale revision of mathematics, of which the revision of logic
consisting in the rejection of excluded middle and double negation
elimination was just a small part. Its mathematical pedigree, however,
ensured continued interest in intuitionistic logic by mathematical
logicians throughout the 20th century. Although it was long considered
a somewhat niche area of mathematics, it has recently become of
central importance. This is due to several factors. One is the
fact that higher-order versions of intuitionistic logic are
sufficiently expressive to develop large parts of mathematics. Another
is that computer-assisted proof systems have matured to the point
where they can be and are being used by mainstream mathematicians to
formalize mathematical proofs. Many of these proof systems use
versions of intuitionistic type theory.

Between the early interest on the part of mathematics in formal logic
in the 1920s and 1930s and the current renewed interest in
computer-aided mathematical proof systems, mathematical logic was not
seen as exactly central to mainstream mathematics. This stands in
stark contrast to the situation in computer science. The theory of
computation itself is an outgrowth of the early advances in
meta-mathematics. Hilbert and his students pursued two main goals in
their investigation of formal logic in the 1920s. The more well-known
is the aim of what's called ``Hilbert's program'': to find elementary
(``finitary'') consistency proofs of axiomatized systems of
mathematics \citep{Zach2023a}. The less well known is the decision
problem: to find an algorithm that decides if a given formula of
predicate logic is a theorem. The negative solution to this problem,
given almost contemporaneously by \citet{Church1936} and
\citet{Turing1937}, marks the beginning of computability theory. The
study of models of computability in theoretical computer science, and
even the study of the complexity of algorithms, is continuous with
this development.

As important as the development of models of computability and
computational complexity is, this development has been for the most
part unrelated, both conceptually and historically, to the development
and philosophy of pure logic. However, philosophical logic has made a
different and perhaps more significant impact in computer science
other than by giving birth to computability theory. Formal logical
systems, and the methods of proof and model theory developed for them,
are used all over the place in computer science. Logical languages,
their proof systems, and semantic frameworks for them have numerous
applications, from theoretical to industrial. Their use is ubiquitous
in computer science. \citet{Halpern2001} have called this the ``unusual
effectiveness of logic in computer science.'' 

Perhaps the earliest and simplest example is the use of Boolean
algebras in the theory of switching circuits \citep{Shannon1938}. But
it was also philosophical logicians who have had significant and
lasting impact. Some other early examples are Quine's work on circuit
minimization (the Quine-McCluskey algorithm), Putnam's work in
automated theorem proving (the DPLL algorithm), and the work of proof
theorists like Prawitz and Martin-Löf on intuitionistic type theories
which now underlies typed programming languages. Systems of modal
logic (especially temporal, epistemic, and deontic systems) are used
for knowledge representation, specification of circuits and programs,
planning in AI and robotics. All of this builds on the pioneering work
of philosophers who developed modal logics and their model and proof
theory, such as Kripke, Hintikka, Prior, Stalnaker, Lewis, von Wright,
Segerberg, Fine, van Benthem, to name just a few. Many-valued logics
find applications, inter alia, in program semantics and reasoning with
imprecise information (``fuzzy logic''). In what follows, we survey
some of the most significant such contributions of logic to
mathematics and computer science.

\section{Logic or mathematics?}

Mathematical logic is widely considered a subfield of mathematics, and
with good reason. It and its traditional subdisciplines (set theory,
model theory, proof theory, and computability theory) can all be found
in the Mathematics Subject Classification, for instance, and papers in
any of them appear in mainstream mathematics journals. As such, one
may wonder where to draw the boundary between logic, as
understood by philosophers, and logic as a mathematical discipline.
This question cannot, and also need not, be settled in a survey paper.
But some preliminary delineations can be made, if only to forestall
complaints about the scope of the discussion below.

Formal logic can be thought of as the study of certain formal
languages, their semantics and proof theory. It is this aspect of
logic which underlies its power and usefulness in other fields,
including, but not only, mathematics. In its modern form, this is
where logic comes from. Its syntax, semantics, and proof theory were
developed primarily to deal with mathematical theories; formal logic
arose as the general theory of axiomatic systems. In order to be
successful, the language has to be complex and expressive enough to
describe its target, and the relationship between its semantics and its
proof theory had to be made clear. In the first instance, from the
philosophical logician's perspective, this amounts to the question
about the relationship between consequence and provability. If all
goes well, they coincide: the proof theory of the logic is both sound
and complete for its semantics. 

\emph{Proofs} of results such as the soundness and completeness
theorems require mathematical methods that go beyond logic. For
soundness, we at least need mathematical induction (on the length of
proofs). For completeness, already in the case of first-order logic,
we need at least a weak form of the axiom of
choice.\footnote{Specifically, for a countable language, what's
required is the weak K\"onig lemma (every infinite binary tree has an
infinite branch). In Henkin-style proofs, it's needed for Lindenbaum's
Lemma.} Not a lot can be established \emph{about} logic without making
use of non-logical, mathematical principles and methods. The case is
similar in other areas of philosophical logic, e.g., the theory of
definitions or of paradoxes. Quite often questions can only be made
precise and solutions provided once we switch to a mathematical
perspective, e.g., that of model theory, proof theory, arithmetic, or
set theory.

Logic can also be thought of as the study of certain logical objects
and properties, such as concepts, propositions, identity, and truth.
Attempts to provide a \emph{logical} foundation for mathematics took
logic in this sense. What Frege and, following him, Whitehead and
Russell attempted to do was to show that mathematics (at least
arithmetic) can be reconstructed in a formal system with only logical
primitives and only assuming logical axioms. Their work showed that
one can accomplish a lot as far as the first desideratum is concerned.
The system of Frege's \emph{Grundgesetze} and Whitehead and Russell's
\emph{Principia Mathematica} only use logical primitives: logical
connectives and quantifiers, with quantifiers ranging over concepts or
propositional functions, respectively. However, Frege's system used an
axiom which was famously inconsistent, and the system of
\emph{Principia} used axioms (reducibility, choice, and infinity)
whose logical status is in doubt.

We see that there is no clear dividing line between logic and
mathematics. Once we move from purely metaphysical questions in the
philosophy of logic to questions about the properties of this or that
logic, we must make use of mathematics---even if the target of study
is a logical system or problem far removed from mathematics. In what
follows, we put emphasis on applications of logic in mathematics and
computer science that concern questions and results which are made
possible, easier, or more informative by the use of logical methods.

Set theory and computability theory will be given somewhat short
shrift. One might well consider these as subfields of logic (in the
philosopher's sense) and not just as areas of mathematics and computer
science with closer historical and conceptual ties to logic than
others. Their omission here does not deny that they are of central
importance to mathematics and computer science as well as to the
philosophy of mathematics and computability. There are of course also
many of examples of questions and results in set theory and
computability theory, understood as mathematical disciplines, that
make essential use of logical methods.

\section{Logic and mathematics}

\subsection{Proof theory}

In the wake of the failure of the Frege/Whitehead/Russell logicist
project to provide a foundation of classical mathematics on purely
logical grounds (see below), and faced with challenges to classical
mathematics from critics such as Brouwer, Poincaré, and Weyl, David
Hilbert followed a different strategy. This strategy pursued two aims.
The first was a formalization program. Hilbert proposed to formalize
classical mathematics in a system similar to that proposed by
Whitehead and Russell, essentially, a first-order theory with a
classical proof system. In fact, in the course of doing so, Hilbert
isolated first-order classical logic for the first time. In contrast
to the logicists, however, he allowed non-logical primitives in the
language. This included constants, functions, and predicates for
mathematical objects, properties, and relations. E.g., a system
suitable for number theory could include $0$, $1$, $+$, $<$ as
primitives. Where the logicists sought to prove the induction
principle from logical axioms, Hilbert took it as a postulate of the
system, which needed no proof.\footnote{\citet{MancosuZachBadesa2009}
give a historical overview of the development of mathematical logic
from \emph{Principia} to the 1930s.}

The second aim was that of finding consistency proofs for the systems
so established, i.e., to secure classical mathematics by showing
that, contrary to what its critics had claimed, classical mathematics
was not self-contradictory. Of course, it could only persuade those
critics if the methods used to establish consistency were themselves
methods the critics accepted. Thus, Hilbert hoped for not just any
proof of consistency, but a ``finitary'' proof of consistency. No
exact definition of ``finitary method'' was given. However, it was
clear that the proofs should avoid any mention of infinite objects,
e.g., only mention finite sequences of symbols. They should also be
constructive, e.g., use computable transformations of such sequences.
A standard pattern of such a proof describes (1)~a finitary
transformation of any proof in the formal system of classical
mathematics to one in a simple form, (2)~a finitary proof that no
proof in simple form can be a proof of a contradiction. Consistency
would follow, since if classical mathematics were inconsistent, there
would be a formal proof of a contradiction in it. The procedure given
in (1)~would transform it into a simple proof of a contradiction.
(2)~establishes that this is impossible, on finitary grounds.

It is widely (though not universally) accepted that Hilbert's program
cannot be carried out. The reason (and proof) was provided by Gödel.
Gödel showed that the kinds of operations on sequences of symbols
required for a consistency proof can be simulated in arithmetic, i.e.,
we can associate natural numbers with sequences of symbols in such a
way that properties of and operations on such sequences can be
captured in the formal system of arithmetic. One such property is that
of being a correct proof of a formula (in a system~$T$ of classical
mathematics). It corresponds to an arithmetical predicate
$\mathit{Prf}(x,y)$. $\mathit{Prf}(n,m)$ holds (and is provable
in~$T$) whenever $n$ and $m$ are numerical codes of a proof and the
formula it proves, respectively. Since $T$ is a formal system that
captures all of classical mathematics, it should also prove whatever
can be proved classically (which includes the finitary methods). In particular,
if we can prove that there is no proof of~$\bot$ (a contradiction)
in~$T$, then we should also be able to prove that no number is the
code of a proof of~$\bot$, and so $T \vdash \lnot \exists
x\,\mathit{Prf}(x, \bot)$. Gödel's second incompleteness theorem
states that if $T$ is a formalized theory (the property of being a
correct proof is decidable), is consistent, and includes elementary
arithmetic (and hence can formalize Gödel ``coding''), then $T \nvdash
\lnot \exists x\,\mathit{Prf}(x, \bot)$.\footnote{The second
incompleteness theorem also depends on the provability predicate
satisfying certain ``provability conditions.''}

Even though Hilbert's program was unsuccessful in reaching its main
aim to establish the consistency of classical mathematics, it was
overall extremely fruitful. The first major advance was the
identification of first-order logic as a suitable system for the
formalization of mathematics. This provided a framework for almost all
future work in mathematical foundations. Axiom systems for set
theories (Zermelo-Fraenkel, but also others like von
Neumann--Bernays--Gödel) are formulated in first-order logic. Axiom
systems for theories investigated by algebraists and geometers are
first-order. Even theories of analysis, a.k.a., ``second-order
arithmetic'' are (two-sorted) first-order systems. This paved the way
for the study of large classes of mathematical theories using the
methods of model theory, proof theory, and set theory.

The second payoff of Hilbert's program was the development of proof
systems. Hilbert himself used axiomatic systems along the lines of
those used by Frege and Russell: a few inference rules, such as modus
ponens, plus axioms for the propositional connectives and quantifiers.
Hilbert's own contribution to the development of logic, aside from the
focus on first-order systems, was the introduction of the epsilon
calculus \citep{AvigadZach2002}. It fell to others at Göttingen and
elsewhere to develop systems more amenable to proof theoretic study:
natural deduction \citep{Gentzen1934,Jaskowski1934} and the sequent
calculus \citep{Gentzen1934}, leading also to tableaux
\citep{Beth1955a}.

Gödel's theorems made it clear that consistency proofs for
mathematical theories will usually require principles that themselves
cannot be proved in the theory. These principles typically take the
form of principles of induction along certain computable well-orders.
Subsequent work in proof theory of mathematical theories thus yields a
classification of the strength of various theories according to their
so-called proof-theoretic ordinals. The first to provide such a
classification was Gentzen. He gave a consistency proof of first-order
arithmetic $\mathbf{PA}$. It uses an idea similar to the one he used
in his proof that the cut rule can be eliminated from proofs in the
sequent calculus. The result is a transformation of proofs in
$\mathbf{PA}$ to simple proofs that avoid cuts on non-atomic formulas
and applications of the induction rule. The proof that this
transformation always terminates requires induction up
to~$\varepsilon_0$.\footnote{The ordinal $\varepsilon_0$ is the limit
of the ordinals $\omega$, $\omega^\omega$, $\omega^{\omega^\omega}$,
\dots. No infinite ordinals are needed in the proof, but it uses
``ordinal notations,'' and these sequences are ordered by a
$<$-relation which is isomorphic to the ordering of ordinals $<
\varepsilon_0$.} To this day, mathematical proof theorists investigate
the proof theoretic strength of stronger and stronger system following
the same pattern.\footnote{See \citet{MancosuGalvanZach2021} for an
accessible introduction to Gentzen's consistency proof,
\citet{RathjenSieg2023} for a recent survey, and \citet{Arai2020} and
\citet{Pohlers2009} for technical introductions to ordinal analysis.}

Ordinal analysis is a kind of reduction of one system to another. A
consistency proof describes a transformation of proofs in one theory,
$T$, to another theory,~$T'$. In the case of Gentzen's analysis of
$\mathbf{PA}$, the transformation is of proofs in first-order
arithmetic to proofs in a simple theory~$T'$. To establish
consistency, the transformation only has to work for (hypothetical)
proofs of a contradiction. In practice, however, the transformations
can be made to work for more complex statements. E.g., an extension of
Gentzen's procedure also works for theorems of the form $\forall
x\exists y\,A(x,y)$ with $A(x,y)$ not containing any unbounded
quantifiers. We can think of formulas of this form as expressing that
the function~$f$, defined by $f(x) = {}$ the least $y$ such that
$A(x,y)$, is total. If $\mathbf{PA} \vdash \forall x\exists
y\,A(x,y)$, then the proof transformation applied to this proof yields
a computation of the function~$f$ using \emph{recursion}
along~$\varepsilon_0$. An ordinal analysis of a theory~$T$ thus
provides two kinds of information. First, it identifies an
ordinal~$\alpha$ so that induction along~$\alpha$ suffices to
establish the consistency of~$T$. Second, it characterizes the
computable functions that $T$ proves to be total (the ``provably
computable functions'' of~$T$) as the $\alpha$-recursive
ones.\footnote{See \citet{RathjenSieg2023}, Appendix~F, for an
explanation of provably computable functions and
\citet{SchwichtenbergWainer2011} for a textbook treatment. The idea
goes back to \citet{Kreisel1952a}, who proved that the provably
computable functions of~$\mathbf{PA}$ are the
$\varepsilon_0$-recursive ones. Mints, Parikh, and Parsons made early
contributions to this line of research, which also led to the study of
theories of bounded arithmetic connecting proof theory and
computational complexity theory; see \citet{Buss1998a}.}

As this example shows, (extensions of) consistency proofs can provide
additional, mathematically useful, information. Suppose that we have a
proof of $\forall x\exists y\,A(x,y)$. What more do we know than that
this is simply true, if we know that this proof can be carried out
in~$\mathbf{PA}$? By the above result, we know that for any~$x$, the
witness~$y$ is bounded by an $\varepsilon_0$-recursive function. In
general, proof theoretic methods can and have been used profitably to
extract information from proofs that these proofs do not obviously
contain. In fact, the proof of $\forall x\exists y\,A(x,y)$ may be
non-constructive. We can make it constructive if we know that it can
be carried out for a system for which suitable proof theoretic methods
are available. This idea is called ``proof mining'': for proofs of
theorems of a certain restricted form (e.g., $\forall x\exists
y\,A(x,y)$), a proof in a restricted theory (e.g., $\mathbf{PA}$)
yields a constructive bound on $y$ depending on~$x$ (e.g., given by an
$\varepsilon_0$-recursive~$f(x)$). The idea originally goes back to
Kreisel. One of the earliest applications was Luckhardt's
\citeyearpar{Luckhardt1989} improvement of bounds in Roth's theorem on
the number of rational approximations to irrational algebraic numbers.
The methods used in proof mining are usually based not on Gentzen-type
consistency proofs, but on consistency proofs using Herbrand
expansions, realizability, or functional interpretations (i.e.,
versions of Gödel's \citeyearpar{Godel1958} \emph{Dialectica}
interpretation of arithmetic extended to stronger
systems).\footnote{See \citet{Kohlenbach2008} for a textbook treatment
and \citet{Kohlenbach2023} for a recent survey.}

Transformations of proofs of theorems of a certain complexity in one
system into proofs in another system are not only used in ordinal
analysis and proof mining. Proof theoretic methods have also helped to
clarify the relationship between various foundational systems. By
foundational system we mean a mathematical theory which is directly
justified by a certain position in the philosophy of mathematics.
Proof theorists and philosophers of mathematics have identified a
number of such systems, all subsystems of second-order arithmetic.
First-order Peano arithmetic itself, e.g, is justified if one accepts
countable infinities. But second-order arithmetic is not so justified,
because it quantifies over sets of natural numbers, and there are
uncountably many of those. Yet, if we suitably restrict the set
existence assumptions, and restrict ourselves to theorems of a certain
complexity, then proof theoretic reductions yield the conservativity
of the stronger theory over the weaker for the theorems so restricted.
E.g., $\mathbf{ACA}_0$ is the theory in which the comprehension axiom
schema
\[\exists X\forall y(y \in X \leftrightarrow A(y))\]
is restricted to arithmetical~$A$ (i.e., not containing
quantifiers over sets) and induction is assumed in the form
\[\forall X[(0 \in X \land \forall y(y \in X \rightarrow y+1 \in X))
\rightarrow \forall y(y \in X)].\] It is not justified on countably
infinitary grounds. Yet, there is a transformation of proofs of
arithmetical theorems (i.e., those not quantifiying over sets) in
$\mathbf{ACA}_0$ to proofs in~$\mathbf{PA}$. Similar reductions are
known for infinitary systems to finitary systems, non-constructive
systems to constructive systems, and impredicative systems to
predicative systems (see
\citealt{Feferman1988a,Feferman1992,Feferman1993a}). Some theories
studied in this context also play a role in the program of reverse
mathematics. Very much like the very first applications of logic to
geometry, reverse mathematics attempts to characterize the minimal
mathematical assumptions required to prove certain theorems of
classical analysis. Typically, this involves identifying a theory $T$
which proves the theorem, and showing that the theorem, together with
a weak base theory, proves the axioms of~$T$. The weak base theory in
this case must be sufficiently strong to formulate analysis, e.g., to
speak about sequences of real numbers. As an example, the
Bolzano--Weierstrass theorem (every bounded sequence has a convergent
subsequence) is equivalent to the theory~$\mathbf{ACA}_0$ over the
base theory~$\mathbf{RCA}_0$.\footnote{See \citet{Eastaugh2024} for a
survey of reverse mathematics, \citet{Simpson2009} for a comprehensive
introduction to subsystems of second-order arithmetic, and
\citet{DeanWalsh2017} for a history of their development.}

\subsection{Model theory}

The most fundamental results about any logic concern the relationship
between semantic entailment, $\Gamma \vDash A$, and provability
$\Gamma \vdash A$. The soundness theorem states that provability
implies entailment; the completeness theorem that entailment implies
provability. These results are also of fundamental importance for the
application of logic to axiomatic theories, especially in mathematics.
From its beginnings, logic was used to codify the primitives used in a
mathematical theory, set down some fundamental truths about these
primitives as axioms, and use formal proof to derive consequences from
these axioms. But mathematics is not just interested in what can be
derived from the axioms, but also what structures realize the axioms,
or whether there indeed are structures that realize them. It was
reasonably clear to mathematicians before \citet{Tarski1936b} that
$\Gamma \vDash A$ must mean that $A$ is true in every realization
(``model'') of the axioms in~$\Gamma$, even if they did not have a
precise definition.

With the insight that a vast array of interesting mathematical
structures can be described using the formal language of first-order
logic, studying realizations of such sets of sentences turned out to
be just a different way of studying the structures. The soundness and
completeness theorems in the form above guarantee that provability
from the axioms exactly captures truth in all structures realizing the
axioms. But as anyone who has seen an actual proof of the completeness
theorem knows, it also provides something else: a guarantee that
structures exist that realize any consistent set of axioms. Proofs of
completeness are almost always model existence theorems of the form:
if $\Gamma$ is consistent (i.e., $\Gamma \nvdash \bot$) then $\Gamma$
has a model~$\mathbf{M} \vDash \Gamma$. In a sense, the completeness
theorem justifies Hilbert's conviction that consistency is all that
counts in mathematics: as long as a theory is consistent, it is a
legitimate object of mathematical study.

In some cases, mathematicians are interested in specific structures:
geometry is interested in Euclidean space, number theory in the
natural numbers, analysis in the reals and complex numbers. In many
other cases, mathematicians are interested in large classes of
structures: algebraists are not usually interested in a single field,
but in all fields, or all finite fields of a certain characteristic,
or all algebraically closed fields. Topologists aren't interested in a
single space, but all locally compact spaces, or all Hausdorff spaces.
Logic, via the results and tools of model theory, applies to these two
kinds of cases differently.

In the first case, when there is a single intended structure, results
from logic are limitative. First-order languages are not able to
describe \emph{only} a single infinite structure (such as the natural
numbers or the real number field), not even up to isomorphism. The
relevant result is a corollary of the completeness theorem: the
compactness theorem. In the version most familiar to philosophical
logicians, it is the property of entailment that $\Gamma$ entails $A$
only if already some finite subset $\Gamma_0 \subseteq \Gamma$
entails~$A$. Suppose that $\Gamma$ entails $A$. By completeness,
$\Gamma$ proves $A$. But any proof can only make use of finitely many
sentences~$\Gamma_0$ in~$\Gamma$, and by soundness, $\Gamma_0 \vDash
A$. For its application in mathematics, the model-theoretic version is
more relevant: if every finite subset of $\Gamma$ is satisfiable, then
$\Gamma$ itself is satisfiable. Here's how this shows that, e.g., the
natural numbers cannot be characterized in a first-order language:
Take any set of sentences~$\Delta$ true in~$\mathbb{N}$ (e.g., the set
of \emph{all} sentences true in~$\mathbb{N}$). Let $\Lambda$ be the
set $\{c \neq \overline{n} : n \in \mathbb{N}\}$, where $c$ is a new
constant symbol, and $\overline{n}$ is a term naming the number~$n$,
e.g., $0+1+\cdots+1$ with $n$ $1$'s. Every finite subset of $\Delta
\cup \Lambda$ is satisfiable: pick a large enough number $k$ for the
constant~$c$. But in a model of $\Delta \cup \Lambda$, the referent of
$c$ must be a number different from every ``natural'' number. Such
non-standard models of arithmetic cannot be isomorphic
to~$\mathbb{N}$.

A similar proof applies to any first-order theory with infinite
models. It can be used to prove the upward Löwenheim--Skolem--Tarski
theorem that every theory with countably infinite models has models of
any infinite cardinality. It can be fruitfully applied to
other theories of mathematical interest, such as the theory of the
real numbers. The real numbers do not contain infinite or
infinitesimal numbers, i.e., numbers~$r$ such that $\left|r\right| >
n$ or $0 < \left|r\right| < 1/n$ for all natural numbers~$n > 0$.
Infinitesimals, i.e., infinitely small but non-zero quantities, were,
however, used for a long time in the development of the calculus.
Berkeley famously ridiculed them as the ``ghosts of departed
quantities.'' The arithmetization of the calculus in the 19th century
replaced infinitesimals by the
$\varepsilon$-$\delta$ definition of limits. Using methods of model
theory, \citet{Robinson1966} was able to provide a rigorous
development of non-standard analysis and show that it is possible to
work consistently with infinitesimals. (It is easy to show the
existence of non-Archimedean fields that make the same statements true
as the ordinary real numbers using a compactness argument. Robinson
used the more sophisticated method of ultraproducts.) Since then,
non-standard analysis has had a following in mathematical circles for
pedagogical reasons. Methods of non-standard analysis have also been
fruitfully applied to develop new mathematical theories and prove new
results.

Results such as Robinson's, as well as most applications of model
theory to algebra, require mathematical, i.e., not purely logical
methods or results. Nevertheless, some important results in
mathematics could not have been achieved without bringing
formalization and logic into play. We'll give two more examples, both
establishing the decidability of important algebraic theories. The
first is the theory of real closed fields. (An ordered field is real
closed if it contains all positive square roots and roots for every
polynomial of odd degree.) The theory of real closed fields is the set
of sentences true in all real closed fields, or, equivalently, in the
real numbers. That this theory is decidable is a surprising fact:
after all, by the theorems of Church and Turing, the theory of the
natural numbers is not decidable. The result implies that elementary
geometry is decidable (by reducing it, via Cartesian coordinates, to
statements about real numbers). It was established by
\citet{Tarski1948}, using the method of quantifier elimination. It is
this step that requires the use of logic: we have to formalize the
language and then find a way to show that every formula is equivalent
to a quantifier-free formula. A decision method for the quantifier
free formulas then provides a decision method for the entire language.
It, however, also requires mathematics to show that the equivalence
holds (in Tarski's case, a generalization of Sturm's theorem about the
real roots of polynomials). 

A second example is the Łos--Vaught test \citep{Los1954,Vaught1954}.
It states that if a theory is categorical in some infinite
cardinality~$\kappa$ and has no finite models, then it is complete. (A
theory $T$ is complete if, for any $A$, either $T \vdash A$ or $T
\vdash \lnot A$. It is categorical in cardinality $\kappa$ if any two
models of cardinality $\kappa$ are isomorphic.) The proof is simple,
and uses only logical methods: Suppose $T$ is not complete, i.e.,
there is some $A$ such that $T_1 = T \cup \{A\}$ and $T_2 = T \cup
\{\lnot A\}$ are consistent. Since $T_1$ and $T_2$ are consistent,
they each have models. Since $T$ has no finite models, they both have
infinite models. By the upward Löwenheim--Skolem--Tarski theorem, they
have models of cardinality~$\kappa$. But $T$ only has a single model
(up to isomorphism) of cardinality~$\kappa$, and no model can make
both $A$ and $\lnot A$ true. If $T$ is an axiomatizable complete
theory, it has a decision procedure: enumerate all proofs from~$T$
until we find either a proof of $A$ or a proof of $\lnot A$. By the
completeness of $T$, this procedure eventually terminates. Armed with
this test, it is easy to prove completeness and hence decidability of
a theory if it is a (mathematical) fact that there is a unique
structure of a given infinite cardinality. For instance, by a theorem
of Cantor, every countable densely ordered set without endpoints is
isomorphic to the rational numbers. Hence, the theory of dense linear
orders without endpoints is $\aleph_0$-categorical, hence complete,
hence decidable. Similarly, by a theorem of Steinitz, there is (up to
isomorphism) only one algebraically closed field of characteristic~$p$
and any uncountable cardinality~$\kappa$. Hence, the theory of
algebraically closed fields of characteristic~$p$ is decidable.

Any two models of a complete theory are elementarily equivalent, i.e.,
they make the same first-order sentences true. The Łos--Vaught test
thus also shows that, e.g, any two algebraically closed fields of
characteristic~$p$ are elementarily equivalent. An important feature
of basic model theory is that elementary equivalence of two structures
does not imply that they are isomorphic, even if they are of the same
cardinality. Mathematicians (or, model theorists aiming to obtain
results in algebra) can exploit this fact to infer the truth of a
first-order statement in one structure from its truth in another,
elementarily equivalent structure. This is yet another example of how
mathematical results can be obtained, and sometimes only obtained, by
formalizing the theories involved in a logical language and applying
general principles about such languages and their models. For
instance, the elementary equivalence of structures was used
essentially in the proof of the Ax--Kochen theorem, which disproved a
long-standing algebraic conjecture of Artin
\citep{AxKochen1965}.\footnote{See \citet{Hodges2023} for an overview
of model theory, and \citet{Hodges1993} and \citet{Marker2002} for
comprehensive, technical introductions. See \citet{ButtonWalsh2018}
for applications of model theory in philosophy.}

\subsection{Logic and mathematical foundations}\label{formal-math}

Two of the most ambitious projects to provide a foundation for
mathematics were those of \citet{Frege1893} and
\citet{WhiteheadRussell1910a}. Their aim was to establish that the
axioms of arithmetic (and in Whitehead and Russell's case, a good deal
more) could be proved from purely logical principles. That is, they
avoided straightforwardly non-logical, mathematical primitives and
axioms, e.g., numbers or sets. Frege attempted to carry out the
reduction of arithmetic to logic in what we would now call higher-order logic. In addition to quantification
over objects, Frege's system allowed quantification over concepts,
concepts that themselves apply to concepts, etc. He successfully
provided a definition of natural number in \citep{Frege1884}. In very
broad strokes, this definition made use of three ideas. The first is
what's now called Hume's law: two concepts have the same number of
objects falling under them if, and only if, they can be put into a
one-to-one correspondence (i.e., are ``equinumerous''). The second is
what's called abstraction. Given a criterion of sameness in some
respect~$X$, such as that provided by Hume's law for sameness of
cardinality, we can introduce higher-type concepts that apply to all
and only those concepts that are the same in respect~$X$. Number
concepts then are concepts that apply to both $F$ and $G$ iff $F$ and $G$
are equinumerous. The concept belonging to~0 is the concept ``is
equinumerous with the concept that has nothing falling under it.''
Numbers as higher-type concepts can then be related to one another in
certain, logically definable ways, e.g., the concept $n$ is related to
the concept $n+1$ in the successor relation. The concept ``is a
natural number'' then is the domain of the transitive closure of the
successor relation, starting from~$0$. That is the third idea: the
logical definition of the ``ancestral'' of a relation. In
\emph{Grundgesetze} \citeyearpar{Frege1893}, Frege additionally
introduced extensions of concepts, which themselves are objects of the
lowest type. It is this additional requirement which rendered his
theory susceptible to Russell's paradox, and hence inconsistent. The
neo-logicist tradition in the philosophy of mathematics has worked on
turning Frege's ideas into a consistent system. Neither Frege's work,
nor that of the neo-logicists, however, has made an impact in
mathematics.\footnote{See \citet{Cook2023} and \citet{Zalta2023} on
Frege's foundations and logical system, \citet{Hale2001} on the
neo-logicist project, and \citet{BoccuniSereni2024} in this volume on the
philosophical import of abstraction principles.}

The contradiction in Frege's system was famously found by Bertrand
Russell in 1901, who, together with Alfred North Whitehead, proceeded
to develop Frege's ideas consistently in \emph{Principia mathematica}
\citeyearpar{WhiteheadRussell1910a}. To this end, they developed a
system of logic, the ramified theory of types, in which propositional
functions (predicates) can themselves be used as arguments and
quantified over. Starting from a basic type of objects, we get
propositional functions of objects, propositional functions that apply
to both objects and propositional functions that apply to objects, and
so on. However, a propositional function can only be applied to one of
lower type. This restriction eliminates the possibility of formulating
Russell's paradox. Even though the notation in \emph{Principia} is
littered with set-theoretic symbols like $\in$ and $\{ x \mid \dots
\}$, it does not take sets or classes as primitive: these symbols are
defined using contextual definitions that only involve propositional
functions and logical vocabulary. (E.g., ``$x \in \{ y \mid A(y)\}$''
is reduced to $A(x)$.) This ``no-class theory'' is a central feature
of \emph{Principia}, which arguably makes it a purely logical system
and not a theory of sets or classes. Although the three volumes of
\emph{Principa} succeed in developing not only the theory of natural
numbers, but even the theory of ordinal and cardinal arithmetic as
well as elements of analysis, it ultimately failed to deliver a
\emph{logical} foundation. The reason is that it made use of two
axioms that cannot be accepted as purely logical, namely the
multiplicative axiom (a version of the axiom of choice) and an axiom
of infinity.

A drawback of the original system was also its system of ramification
of propositional functions: in determining the type of a propositional
function, it not only took into account the types of its arguments,
but also the types of propositional functions that were quantified
over in its definition. The axiom of reducibility, which states that
every propositional function is coextensive with a predicative
propositional function, essentially undoes the ramification. Following
a suggestion of Chwistek and Ramsey, Russell proposed to develop the
system of \emph{Principia} simply on the basis of an un-ramified
(``simple'') theory of types. Simple theories of types were first
described by \citet{Carnap1929} and \citet{Church1940}. Although
simply typed logics were studied by logicians subsequently, they also
did not make an impact on mainstream mathematics for a long time: by
the mid-20th century and up to now, the preferred foundational
framework, as far as mathematicians were concerned, was
Zermelo--Fraenkel set theory.\footnote{See \citet{LinskyIrvine2024}
for an overview of \emph{Principia} and \citet{Coquand2022} for a
survey of type theories.}

Type theories \emph{have} made a significant impact on mathematics
more recently, however. \citet{Milner1972} described a type theory
related to Church's simple type theory: the LCF system (logic for
computable functions), based on an earlier system described by Scott
in 1968 \citep{Scott1993}. Another crucial contribution to this
development was the work of Martin-Löf on intuitionistic type theories
\citep{Martin-Lof1975,Martin-Lof1982,Martin-Lof1984}. It is a
higher-order version of intuitionistic logic which essentially uses
the Curry--Howard correspondence between propositions and types and
between proofs and programs (see below). An impredicative system
similar to Martin-Löf's and incorporating features of Girard's
system~$F$ \citep{Girard1971} is the calculus of constructions
\citep{CoquandHuet1988}. These systems form the basis of computerized
proof assistants: LCF is the basis of HOL and Isabelle/HOL, Martin-Löf
type theory is the basis of Nuprl and Agda, and the calculus of
constructions that of Coq (soon to be renamed Rocq) and Lean. (Another
system to mention here is Mizar, although it is based on set
theory.)\footnote{On HOL (\url{https://hol-theorem-prover.org/}) see
\citet{Gordon2000}; on Isabelle (\url{https://isabelle.in.tum.de/})
see \citet{Paulson1989} and \citet{PaulsonNipkowWenzel2019}. Nuprl
(\url{https://nuprl-web.cs.cornell.edu/} was originally developed by
\citet{Constableothers1986}, Agda
(\url{https://wiki.portal.chalmers.se/agda/}) by \citet{Norell2007} on
the basis of the version of type theory given by \citet{Luo1994}, and
Coq (\url{https://coq.inria.fr/}) by Coquand and Huet. On Lean
(\url{https://lean-lang.org/}), see
\citet{deMoura2015} and \citet{deMouraUllrich2021}. \textsc{Mizar}
(\url{http://mizar.org/}) is due to \citet{Trybulec1977}. These systems
were all influenced by de Bruijn's \textsc{Automath}
\citep{deBruijn1970}. \citet{HarrisonUrbanWiedijk2014} provide a survey
of the history of interactive theorem proving.}   These proof
assistants have been used to formalize, and formally verify,
research-level mathematical results. E.g., Gonthier verified the four
color and Feit--Thompson theorems in Coq, Hales verified his proof of
the Kepler conjecture in Isabelle/HOL, and Tao the proof of the
polynomial Freiman--Ruzsa conjecture in Lean. The homotopy type theory
and univalent foundations projects are carrying out their work in Coq
as well.

Proof assistants and formal verification have now entered mainstream
mathematics. Notable mathematicians like the the Fields medalists
Gowers, Scholze, Tao, and Voevodsky have advocated their use, and have
themselves used them. Teams of mathematicians are contributing to
growing ``libraries'' of formally verified mathematical concepts and
results that other mathematicians can draw upon. The formalization of
mathematics was shown to be possible by the development of logic and
foundational systems like \emph{Principia} and Zermelo--Fraenkel set
theory a century ago. Formalization of mathematics has now turned from
an in-principle possibility to something mathematicians actually do,
with significant implications for the practice and philosophy of
mathematics.\footnote{See \citet{Avigad2018,Avigad2024} for recent
surveys and discussion of examples and implictions.
\citet{Gonthier2008} describes the formalization of the four-color
theorem and \citet{Hales2017}that of the Kepler conjecture. See
\citep{hottbook,AwodeyPelayoWarren2013,AwodeyCoquand2013}
on univalent foundations and homotopy type theory. \citet{Elkind2022}
provides a useful introduction to the use of proof assistants for
philosophers, including examples of uses in philosophy. Incidentally,
Elkind is in the process of verifying \emph{Principia mathematica} in
Coq, see \url{https://www.principiarewrite.com/}.}

\section{Logic and computer science}

\subsection{Automated reasoning}

Logic is concerned with proof systems for various consequence
relations, and logicians have developed many such systems. In applying
logic to mathematics, we formulate specific theories which are of
interest to mathematicians, and investigate what can be proved from
them. This includes specific consequences---e.g., is this statement a
theorem of the theory?---and also the properties of theory
itself---e.g., is it consistent? complete? decidable? Applications of
logic in computer science are similar in nature, but there are several
important differences. While in mathematics the number of interesting
theories is (relatively) small and varies little over time, in
computer science it is enormous and constantly changing. Every data-
or knowledge base, every specification of a circuit or program is,
fundamentally, a theory; every time a record, fact, or rule is added
or removed, the theory changes. While mathematicians rarely actually
formalize their theories and theorems (except when formally verifying
results; see the previous section), in computer science they are
almost always formally represented (they are typically stored in some
kind of symbolic format, whether recognizable as a formula or some
equivalent, but computationally more efficient, data structure).
Mathematicians usually have candidate theorems in mind, and wonder
whether they are provable in a particular theory (e.g., is the
Bolzano--Weierstrass theorem provable in $\mathbf{ACA}_0$?). By
contrast, in computer science, one is very often not concerned with
specific theorems, but all theorems of a certain form (e.g., a
database query might ask for which values of~$x$ does $\exists
y\,R(x,y)$ follow from the data). Mathematicians are not very (or at
all) concerned with the time required of finding out if a proof
exists; for computer science efficiency is critical.\footnote{Issues
of computational complexity are tied closely to issues in logic and
the philosophy of mathematics; see \citet{Dean2019,Dean2021}.} Computer
scientists need fast algorithms to answer such questions, and they
need to answer literally millions such questions every day.

It is not surprising then that the search for efficient logical proof
methods began more or less as soon as digitial computers became
available. The very first example of a logical theorem prover was the
Logic Theorist \citep{NewellSimon1956}, which implemented heuristic
proof search in the propositional system of \emph{Principia
mathematica}. This is, of course, not a computationally efficient
approach. There are far better approaches. Philosophical logicians
laid much of the groundwork in the early development of automated
theorem proving methods, both by developing the theoretical
foundations and by implementing and testing them on actual computers.
One approach is based on proof search in calculi that are more suited
to it than axiomatic systems, such as \emph{Principia} or Hilbert's,
namely analytic proof systems like tableaux or the sequent calculus.
The first of these was used by Prawitz
\citep{PrawitzPrawitzVoghera1960}, the second by \citet{Wang1960}. 

A second early approach was based on Herbrand's theorem. The theorem
is originally due to \citet{Herbrand1930} and
\citet{HilbertBernays1939}, but the form used was formulated by
\citet{Dreben1952} and \citet{Quine1955a}. Since any formula is valid
iff its negation is unsatisfiable, instead of giving a method that
shows arbitrary formulas are valid, we can give one that shows they
are unsatisfiable. Recall that any formula of first-order logic has an
equivalent prenex form, e.g., $\forall x\, P(x) \land \exists y\,
\lnot P(y)$ is equivalent to $\forall x\exists y (P(x) \land \lnot
P(y))$. Such a formula is satisfiable iff its Skolem form is
satisfiable. In the Skolem form we remove existential quantifiers and
replace the variables they bind by constants or functions depending on
the preceding universally quantified variables, e.g., $\forall x(P(x)
\land \lnot P(f(x)))$. Herbrand's theorem states that such a formula
is unsatisfiable iff a conjunction of instances of its quantifier-free
matrix is propositionally unsatisfiable. In our example, while $P(a)
\land \lnot P(f(a))$ is satisfiable, the conjunction of the two
instances where we replace $x$ once by $a$ and once by $f(a)$, i.e.,
\[(P(a) \land \lnot P(f(a))) \land (P(f(a)) \land
\lnot P(f(f(a))))\] is unsatisfiable. Herbrand's theorem thus reduces
the problem of dealing with first-order formulas to propositional
logic. The remaining problems are those of finding the suitable
substitution instances, and to efficiently show that a propositional
formula is (un)satisfiable.

For the second question, an efficient procedure was given by
\citet{DavisPutnam1960}. It was subsequently refined by
\citet{DavisLogemannLoveland1962}---the ``DPLL'' method is still at the
core of so-called SAT solvers, programs that can effectively determine
satisfiability of very large propositional formulas. These are in
regular use even in industrial applications, thanks to the fact that
most real-life problems in computer science concern finite domains or
structures, and can thus be described without quantifiers. The first
question turned out to be harder to address: simply trying out all
possible instances doesn't work except for almost trivial cases. A
breakthrough was made by \citet{Robinson1965}.\footnote{Robinson was a
philosopher who received an MA from the University of Oregon under
Arthur Pap and a PhD from Princeton advised by Hempel and Putnam.} His resolution calculus
deals with sets of clauses, i.e., sets (interpreted as disjunctions)
of atomic and negated atomic formulas with free variables (interpreted
as universally quantified). Such sets result naturally from matrices
of formulas in Skolem form: simply transform the matrix into
conjunctive normal form, and distribute the initial universal
quantifiers. Applied to our example $\forall x(P(x) \land \lnot
P(f(x)))$, this yields $\forall x\,P(x) \land \forall y\,\lnot
P(f(y)))$ and the clauses $\{P(x)\}$ and $\{\lnot P(f(y))\}$. The
(binary) resolution rule is the following:
\[\AxiomC{$C \cup \{P(t)\}$} \AxiomC{$D \cup \{\lnot P(s)\}$}
\BinaryInfC{$C\sigma \cup D\sigma$} \DisplayProof\] provided there is
a substitution~$\sigma$ (a ``unifier'') such that $P(t)\sigma =
P(s)\sigma$. The clauses in the premise will always have their
variables renamed so that the variables in the two premises are
disjoint.\footnote{The actual rule is more complicated, since it must
account for the possibility that the ``clash'' between $A$ and $\lnot
A$ involves multiple instances of $P(t)$ and $P(s)$ on each side and
arity of $P$ greater than~$1$, and also accommodate ``factoring'' a
clause. Conversion to clause form can be done more efficiently,
resulting in shorter refutations, using methods alternative to
simple-minded prenexation.} A refutation is a proof of the empty set
from the starting set of clauses; a formula is unsatisfiable if its
corresponding clause set has a refutation. Our example has a very
simple resolution refutation:
\[\AxiomC{$\{P(x)\}$} \AxiomC{$\{\lnot P(f(y))\}$}
\BinaryInfC{$\emptyset$} \DisplayProof\]
where the unifier is $\sigma(x) = f(y)$.

Resolution still forms the core of many general-purpose automated
theorem provers such as the E~prover, Prover9, or
Vampire.\footnote{Systems in practical use extend resolution using
various methods, e.g., paramodulation or superposition to deal with
equality, methods from computer algebra to deal with equational
theories, and often try to simultaneously disprove a conjecture using
model building methods. Examples of such provers are Prover9/Mace4
(\url{https://www.cs.unm.edu/~mccune/prover9/}; the successor to
Otter, \citealt{McCune2005}), the E prover (\url{https://eprover.org};
\citealt{SchulzCruanesVukmirovic2019}) and Vampire
(\url{https://vprover.github.io/}; \citealt{KovacsVoronkov2013}).} One
significant legacy of the resolution method is, however, that it
provides the foundation for declarative (``logical'') programming
languages such as Prolog
\citep{ColmerauerKanouiPaseroRoussel1973,ColmerauerRoussel1996,Kowalski1974}.
Here, resolution is restricted to Horn clauses (clauses with at most
one positive formula). Programs consist of facts (just one positive
atomic formula) and rules. ``Running'' a program involves proving a
formula (the ``goal'') from the program using resolution. Crucially,
the goal may contain free variables, and in finding a resolution
proof, the prover generates a substitution that results in a
refutation. Consider the example ``program'' $\forall x(\lnot P(x)
\lor Q(x))$, $P(a)$. (Here $\forall x(\lnot P(x) \lor Q(x)$ is a rule,
which might be written as $Q(x) \leftarrow P(x)$.) We can prove $Q(a)$
from this, but we can also ask, ``For what values of $y$ is $Q(y)$
provable?'' A resolution refutation of the corresponding clauses
$\{\lnot P(x), Q(x)\}$, $\{P(a)\}$, $\{\lnot Q(y)\}$ will produce a
substitution $\sigma(y) = a$ which provides the answer. This principle
underlies not only logic programming, but can also be used for
reasoning in query languages for databases and systems for knowledge
representation and formal ontologies such as description
logics.\footnote{For textbook treatments of the resolution method, see
\citet{Fitting1996a} and \citet{Leitsch1997}. Description logics
\citep{BaaderHorrocksSattler2008} provide a formal framework that is
closely connected to ``real-life'' formalisms for describing
ontologies, i.e., the relationships between objects and concepts in
various domains. Two prominent examples of the latter are SNOMED/CT
for the medical field, and the Web Ontology Language (OWL2;
\url{http://www.w3.org/TR/owl2-overview/}). Although reasoning systems
for description logics are often specific to the language to take
advantage of restrictions in the syntax, they all make use of the same
basic logical foundations (e.g., tableaux or resolution provers). See
\citet{BienvenuLeclereMugnierRoussetEtAl2020} for a survey of
reasoning in formal ontologies.}

The unification principle is crucial for the success of resolution
provers and reasoning in declarative languages. Its usefulness lies in
the fact that the problem of computing unifiers is not only
effectively decidable, but also that the unification algorithm
produces a \emph{most general} solution of which all other solutions
are themselves instances. Unification is necessary not only in
resolution, but in other first- and higher-order proof techniques as
well. It tells us when we can stop because we have reached an axiom.
E.g., consider a search for a proof of the negation of our example
formula $\forall x(P(x) \land \lnot P(f(x)))$ in the sequent
calculus.\footnote{For proof search, a contraction-free calculus is
best suited, such as the system $G3$ of \citet{Kleene1952}.} Applying
introduction rules backwards and leaving witness terms as free
variables we get: \let\fCenter\vdash
\[
  \Axiom$P(y), \lnot P(f(y)), P(x), \forall x(P(x) \land \lnot P(f(x))) \fCenter P(f(x))$
  \RightLabel{$\land$L}
  \UnaryInf$P(y) \land \lnot P(f(y)), P(x), \forall x(P(x) \land \lnot P(f(x))) \fCenter P(f(x))$
  \RightLabel{$\forall$L}
  \UnaryInf$P(x), \forall x(P(x) \land \lnot P(f(x))) \fCenter P(f(x))$
  \RightLabel{$\lnot$L}
  \UnaryInf$P(x), \lnot P(f(x)), \forall x(P(x) \land \lnot P(f(x))) \fCenter$
  \RightLabel{$\land$L}
\UnaryInf$P(x) \land \lnot P(f(x)), \forall x(P(x) \land \lnot P(f(x))) \fCenter$
\RightLabel{$\forall$L}
\UnaryInf$\forall x(P(x) \land \lnot P(f(x))) \fCenter$
\RightLabel{$\lnot$R}
\UnaryInf$\fCenter \lnot\forall x(P(x) \land \lnot P(f(x)))$
\DisplayProof
\]
At this point, computing the unifier $\sigma(y) = f(x)$ tells us that
we can stop: applying the substitution $\sigma$ yields a proof from
the axiom $P(f(x)) \vdash P(f(x))$. This use of unification is
essential in proof assistants, which work by proof search (in suitable
systems of natural deduction). Since those are often higher order
systems, higher-order versions of unification have to be
used.\footnote{The use of unification in proof search was already
suggested by \citet{Prawitz1960} and influenced Robinson, who gave a
first algorithm for it. Higher-order unification is undecidable, as
was shown by \citet{Goldfarb1981}, but more restricted forms are
tractable and often suffice in practice.}

\subsection{Verification of programs and systems}

Since the 1970s, computer scientists have developed many approaches to
and implemented tools for a connected set of problems: verifying that
systems behave the way they are supposed to behave. One version of
this question can be asked about programs: does program $\pi$ compute
the function it is supposed to compute, or more generally: does $\pi$
exhibit the correct input-output behavior? But the question can also
be asked about digital circuits, communication and scheduling
protocols, and a limitless number of other systems that can be
formally described. In order to solve problems of this sort, the
system and its (intended) properties have to be described
(specification), and it has to be proved that the system has the
intended properties or lacks unintended properties (verification).
Methods and ideas of logic are used both in specification (formal
languages) and verification (proof and model theory). Such ``formal
methods'' have developed to the point where they can be used in
large-scale, industrial applications. This is due, on the one hand, to
theoretical advances and the increasing power of modern computers. On
the other hand, systematic testing has turned out to often be
unreliable when applied to complex systems, highlighting the need for
formal verification in addition to testing.\footnote{The two most
famous examples of such failures are the explosion, upon launch, of
the Ariane~5 rocket and the Intel Pentium floating-point bug, each
resulting in hundreds of millions of dollars in losses. Both prompted
the subsequent use of formal methods to verify software and chip
design at Ariane Aerospace and Intel.}

Floyd--Hoare logic \citep{Floyd1967,Hoare1969} is a calculus designed
to prove properties of programs. Its syntax combines a first-order
language to express properties of states, and the syntax of the
programming language used to implement an algorithm. A \emph{Hoare
triple} is an expression $\{A\}\, \pi\,\{B\}$, where $A$ and $B$ are
formulas in the former (the pre- and post-conditions) and $\pi$ is a
program. It states that after program~$\pi$ is run in a state where
$A$ is true, if the program halts then $B$ will be true. The system
has axioms and rules. An example of an axiom is that for variable
assignments,
\[
  \AxiomC{}
  \RightLabel{A}
  \UnaryInfC{$\{A[t/x]\}\, x \mathrel{\mathtt{:=}} t\, \{A\}$}
  \DisplayProof
\]
and examples of rules are
\[
  \AxiomC{$A \to B$} 
  \AxiomC{$\{B\}\, \pi\, \{C\}$}
  \RightLabel{$\to$}
  \BinaryInfC{$\{A\}\, \pi\, \{C\}$}
  \DisplayProof
  \quad
  \AxiomC{$\{A \land C\}\,\pi\,\{A\}$}
  \RightLabel{W}
  \UnaryInfC{$\{A \land C\}\, \mathtt{while}\,C\,\mathtt{do}\,\pi \{A \land \lnot C\}$}
  \DisplayProof
\]
The first allows us to weaken the precondition, and to bring in
background information (from logic and, say, arithmetic). The second
rule says that if the ``invariant'' $A$ is not changed by $\pi$ as
long as $C$ holds, then a \texttt{while} $C$ loop does not change it
either (and that, if and when the loop terminates, $C$ is false).

Suppose we want to verify that the simple program 
\[\mathtt{while}\, x<5\, \mathtt{do}\, x \mathrel{\mathtt{:=}} x+1 \tag{*}\]
computes $x = 5$ if started on $x = 1$, i.e., we want to prove
\[\{x=1\}\, \mathtt{while}\, x<5\, \mathtt{do}\, x \mathrel{\mathtt{:=}} x+1
\,\{x=5\}\]
We have to select a suitable invariant~$A$, i.e., a proposition that
is unchanged by $\pi$ as long as $x<5$ is true. $A$~should also be
implied by our desired precondition $x=1$ and, together
with $x \not< 5$, imply the desired postcondition $x=5$. A suitable
candidate is $x \le 5$. We can prove:
\[
  \AxiomC{$(x \le 5 \land x < 5) \to x+1 \le 5$}
  \AxiomC{}
  \RightLabel{A}
  \UnaryInfC{$\{x+1\le 5\}\,x \mathrel{\mathtt{:=}} x+1\,\{x \le 5\}$}
  \RightLabel{$\to$}
  \BinaryInfC{$\{x\le 5 \land x<5\}\,x \mathrel{\mathtt{:=}} x+1\,\{x \le 5\}$}
  \RightLabel{W}
  \UnaryInfC{$\{x\le 5 \land x<5\}\, \mathtt{while} \,x<5\, \mathtt{do}\, x\, \mathrel{\mathtt{:=}} x+1\, \{x \le 5 \land x\not<5\}$}
  \DisplayProof  
\]
The result (*) follows from the following arithmetical facts:
\begin{align*}
  (x \le 5 \land x < 5) & \to x+1 \le 5\\
  x = 1 & \to (x\le 5 \land x<5)\\
  (x \le 5 \land x \not< 5) & \to x=5
\end{align*}
The first justifies the top left formula in the derivation. The other
two allow us to weaken the precondition to $x = 1$ and strengthen the
postcondition to $x=5$.

This simple example illustrates a number of features of this approach
to program verification. It is a proof-based approach, i.e., there is
no semantic component. As such, it is also formal: we have a formal
language that incorporates logic, but also the syntax of a programming
language as well as expressions for specific domains (in our example,
arithmetic). Proof search is highly indeterminate, e.g., in our
example we had to \emph{guess} the invariant $x \le 5$ for the
\texttt{while} loop. It can also only verify \emph{partial}
correctness: the post condition holds \emph{if} the program
terminates; it does not guarantee termination. For automation, it must
be combined with automated proof methods for the background domains as
well as logic, tools for generating loop invariants, and termination
checkers. (Checking for termination of programs is undecidable in the
general case, but powerful tools that work in many cases exist.)

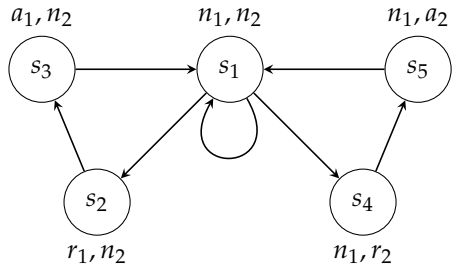
\begin{figure}
\begin{center}
\tikzset{node distance=2.5cm,
every edge/.style={ 
  draw,
  ->,>=stealth,
  auto,
  semithick}
}
\begin{tikzpicture}
  \node[state,label=above:{$n_1, n_2$}] (s1) {$s_1$};
  \node[state, below left of=s1, label=below:{$r_1, n_2$}] (s2) {$s_2$};
  \node[state, left of=s1, label=above:{$a_1, n_2$}] (s3) {$s_3$};
  \node[state, below right of=s1, label=below:{$n_1, r_2$}] (s4) {$s_4$};
  \node[state, right of=s1, label=above:{$n_1, a_2$}] (s5) {$s_5$};
  \draw (s1) edge (s2);
  \draw (s2) edge (s3);
  \draw (s3) edge (s1);
  \draw (s1) edge (s4);
  \draw (s4) edge (s5);
  \draw (s5) edge (s1);
  \draw[->,loop,looseness=7,in=240,out=300] (s1) edge (s1);
\end{tikzpicture}
\end{center}

\caption{First come first served protocol}
\label{protocol}
\end{figure}

Model checking is a different approach to specification and
verification of a wide range of systems. Consider as an example the
following simple-minded protocol for assigning a resource to one of
two agents (an agent might be a program process running on a computer,
a resource might be a device or drive). The agents can be neutral
(doing something else, not involving the resource), they can request
the resource, or they can access the resource. The resource is
assigned ``first come first served,'' i.e., after an agent requests
the resource, they access it, then immediately release it. Only one
agent can request the resource at a time, and it cannot be requested
while it's being accessed.\footnote{Not allowing both agents to even
request the resource simultaneously makes the example very
unrealistic, but we want to keep things extremely simple.} We can use
the propositions $n_i$, $r_i$, $a_i$ for the three possible states
each agent can be in, and the possible states of the system and
transitions between them in the diagram in Figure~\ref{protocol}. The
initial state is $s_1$, where both agents are neutral. One or the
other agents can request the resource, corresponding to states $s_2$
and $s_4$, and once requested, each agent accesses the resource
(states $s_3$ and $s_5$), and then transitions back to a neutral
state~($s_1$).

Anyone familiar with modal logic can recognize this diagram as a
Kripke model. Model checking involves checking whether formulas of
certain modal languages are true in such a model. Actually, it is not
this model itself, but the unravelling of the model starting at the
initial state~$s_1$. The unravelling of the model is a tree with
root~$s_1$ and all other nodes copies of the nodes $s_1$--$s_5$. The
branches of this tree record all the possible (possibly infinite)
paths through the original model, i.e., the possible ways in
which the system can evolve over time.

Naturally the modal logics considered include temporal operators. In
linear temporal logic LTL \citep{Pnueli1977,Pnueli1981}, the basic
operators are $\mathsf{X}\,A$ ($A$ is true in the next state) and $A
\mathrel{\mathsf{U}} B$ ($B$ eventually becomes true, and $A$ is true
until it does). The invention of LTL was influenced by previous work
done by philosophers on logics of time and modal logics generally,
specifically Kripke's semantic framework and the work of
\citet{Prior1967}, \citet{Kamp1968}, and \citet{RescherUrquhart1971}.
Other familiar operators can be defined from $\mathsf{U}$, e.g.,
$\mathsf{F}\,A \equiv \top \mathrel{\mathsf{U}} A$ ($A$ is eventually
true) and $\mathsf{G}\,A \equiv \lnot\mathsf{F}\,\lnot A$ ($A$ is
henceforth always true). These operators are not evaluated with
respect to all possible future states, but only the states along a
branch. It then becomes possible to express properties about the
transition system as modal formulas. For instance, safety properties
say that undesirable states never happen. In our example, we want the
resource never to be accessed by both agents at the same time. This
can be expressed as $\mathsf{G}\, \lnot(a_1 \land a_2)$ (and is true
on every branch). Liveness properties say that something good always,
or always eventually happens. For instance, we want it to be the case
that every agent which requests the resource will be
granted access. We can formalize this as $\mathsf{G}(r_i \to
\mathsf{F}\,a_i)$ (and it is also true on every branch).

Some important properties cannot be expressed in LTL, and this has
lead to extensions of the formalism to include quantification over
branches. Computational tree logic CTL$^*$ \citep{EmersonHalpern1983}
adds the operators $\mathsf{A}\,B$ and $\mathsf{E}\,B$ for this
purpose. Formulas are evaluated relative to branches and
states on them. The formula $\mathsf{A}\,B$ is true at a state~$w$ if
$B$ is true at $w$ relative to all branches through~$w$;
$\mathsf{E}\,B$ if $B$ is true relative to at least one branch
through~$w$. LTL is a subsystem of CTL$^*$, since any LTL formula~$B$
is true at the root iff $\mathsf{A}\,B$ is true at the root
in~CTL$^*$. Here is a fairness property we might want to state (and
verify): it's always true that agent~$2$ \emph{could} access the
resource. We can't formalize that in LTL. The best we can do is
$\mathsf{G}\,\mathsf{F}\, a_2$, which says that it's always true that
agent~$2$ \emph{will} access the resource. But there is a branch in
which only agent~$1$ ever accesses the resource on which this is
false. In CTL$^*$ we can formalize it as
$\mathsf{A}\,\mathsf{G}\,\mathsf{E}\,\mathsf{F}\,a_2$, which is true
for our protocol: for every branch~($\mathsf{A}$), and every future
state on that branch ($\mathsf{G}$), there is a branch through that
state ($\mathsf{E}$) containing a future state ($\mathsf{F}$)
where~$a_2$ is true.

CTL$^*$ contains the system CTL \citep{ClarkeEmerson1982} as a
subsystem (here the branch quantifiers $\mathsf{A}$ and $\mathsf{E}$
can only come immediately before one of the basic modal operators).
CTL is less expressive than CTL$^*$ and does not contain LTL, but is
sufficiently expressive for many applications and has efficient model
checking algorithms. For this reason it is perhaps the model checking
formalism most widely used in practice.\footnote{Temporal logics are
not the only formalisms inspired by philosophers' work on modal logic
that are used in computer science. Others worth mentioning are dynamic
logic \citep{Pratt1976} (closely related to Floyd--Hoare logic),
epistemic logic \citep{FaginHalpern1987,FaginHalpernMosesVardi1995},
and deontic logic \citep{WieringaMeyer1994}. The standard textbook on
model checking, temporal and other logics used in computer science is
\citet{HuthRyan2004}.}

\subsection{Type systems for programming languages}

In the Brouwer--Heyting--Kolmogorov interpretation of
intuitionistic logic, the meanings of a proposition is the set of its
proofs, and proofs are constructions of a certain kind. A construction
(proof) of a conjunction, $A_1 \land A_2$ consists of a pair $\langle
x_1, x_2\rangle$ where $x_1$ is a construction (proof) of $A_1$ and
$x_2$ one of $A_2$. A construction of $A_1 \lor A_2$ is a
construction of one or the other, plus the information which one it is
(i.e., either $\langle 1, x\rangle$ or $\langle 2, x\rangle$). A
construction of $A_1 \to A_2$ is an operation that transforms a
construction of $A_1$ into one of $A_2$. $\bot$ is a proposition that
has no proof (absurdity), and negation $\lnot A$ is defined as $A \to
\bot$.

This interpretation invites a comparison to the natural deduction
rules for intuitionistic logic. We write them in sequent form, to
display the assumptions each formula depends on. The rules for $\land$
are:
\[
  \Axiom$\Gamma \fCenter A_1$
  \Axiom$\Gamma \fCenter A_2$
  \RightLabel{$\land$I}
  \BinaryInf$\Gamma \fCenter A_1 \land A_2$
  \DisplayProof\quad
  \Axiom$\Gamma \fCenter A_1 \land A_2$
  \RightLabel{$\land$E$_1$}
  \UnaryInf$\Gamma \fCenter A_1$
  \DisplayProof\quad
  \Axiom$\Gamma \fCenter A_1 \land A_2$
  \RightLabel{$\land$E$_2$}
  \UnaryInf$\Gamma \fCenter A_2$
  \DisplayProof
\]
We might now add to these rules the information about the
constructions involved, and how the construction itself is built up
from elementary operations. We'll use the notation $t : A$ for ``$t$
is a construction of~$A$.'' The rules then become:
\[
  \Axiom$\Gamma \fCenter t_1 : A_1$
  \Axiom$\Gamma \fCenter t_2: A_2$
  \RightLabel{$\land$I}
  \BinaryInf$\Gamma \fCenter \mathsf{pair}(t_1, t_2) : A_1 \land A_2$
  \DisplayProof\quad
  \Axiom$\Gamma \fCenter t : A_1 \land A_2$
  \RightLabel{$\land$E$_1$}
  \UnaryInf$\Gamma \fCenter \mathsf{proj}_1(t) : A_1$
  \DisplayProof\quad
  \Axiom$\Gamma \fCenter t: A_1 \land A_2$
  \RightLabel{$\land$E$_2$}
  \UnaryInf$\Gamma \fCenter \mathsf{proj}_2(t) : A_2$
  \DisplayProof
\]
The rules now simply record the BHK interpretation: If $t_1$ and $t_2$
are constructions of $A_1$ and $A_2$, respectively, then
$\mathsf{pair}(t_1, t_2)$, the pair-forming operation applied to $t_1$
and $t_2$, is a construction of $A_1 \land A_2$.

For the conditional, the rules become:
\[
  \Axiom$x:A, \Gamma \fCenter t: B$
  \LeftLabel{$x$}\RightLabel{$\to$I}
  \UnaryInf$\Gamma \fCenter \lambda x.t : A \to B$
  \DisplayProof\quad
  \Axiom$\Gamma \fCenter t: A \to B$
  \Axiom$\Gamma \fCenter s: S$
  \RightLabel{$\to$E}
  \BinaryInf$\Gamma \fCenter (ts) : B$
  \DisplayProof
\]
Since in the BHK interpretation, a construction of $A \to B$ is a
procedure for turning a construction of $A$ into one of $B$, the
notation used for the conclusion of $\to$I must be something that
describes a function. That's what the lambda abstract $\lambda x.t$
does: $\lambda x.t$ is a function with argument $x$ that's defined
by~$t$. (Here $t$ in general will contain free occurrences of the
variable $x$.) When we apply a construction of $A \to B$ to a
construction of $A$ we obtain a construction of $B$. Applying a
function~$t$ to $s$ is symbolized as~$(ts)$.

The terms introduced this way are nothing but terms in the lambda
calculus with products. The \emph{typed} lambda calculus is a version
of the lambda calculus where not every well-formed expression is also
well-typed. E.g., $(\mathsf{pair}(x,y)x)$ is not well typed, because
we can only apply functions to arguments, and $\mathsf{pair}(x,y)$ is
not a function, but a pair. The above rules then become rules not for
inferring formulas, but for determining that terms in the lambda
calculus have a certain type. E.g., $\land$I says that if $t_1$ has
type $A_1$ and $t_2$ has type $A_2$, then $\mathsf{pair}(t_1, t_2)$
has type $A_1 \land A_2$.\footnote{In the typed lambda calculus, the
notation for product and sum types is usually $\times$ and $+$ instead
of $\land$ and $\lor$.} The fact that the introduction and elimination
rules of intuitionistic natural deduction are nothing but the typing
rules of the typed lambda calculus is the first part of the
Curry--Howard correspondence.\footnote{It was first observed by Curry
for the conditional rules and combinatory logic. Howard formulated it
for intuitionistic arithmetic and the lambda calculus in a manuscript
from 1969 \citep{Howard1980}. It was independently identified by
\citet{deBruijn1970} and \citet{Reynolds1974}. See
\citet{GirardTaylorLafont1989}, \citet{SorensenUrzyczyn2006}, and
\citet{BarendregtDekkersStatman2013} for in-depth discussions of typed
lambda calculi and the Curry-Howard correspondence.} In it, formulas
correspond to types, and proofs correspond to lambda terms, i.e., to
programs (the lambda calculus is essentially a programming language).
For instance, consider the simple proof of $(A \land B) \to (B \land
A)$:
\[
  \Axiom$x : A \land B \fCenter x : A \land B$
  \RightLabel{$\land$E}
  \UnaryInf$x : A \land B \fCenter \mathsf{proj}_2(x) : B$
  \Axiom$x : A \land B \fCenter x : A \land B$
  \RightLabel{$\land$E}
  \UnaryInf$x : A \land B \fCenter \mathsf{proj}_1(x) : A$
    \RightLabel{$\land$I}
  \BinaryInf$x : A \land B \fCenter \mathsf{pair}(\mathsf{proj}_2(x), \mathsf{proj}_1(x)) : B \land A$
  \RightLabel{$\to$I}
  \UnaryInf$\fCenter\lambda x: (A \land B). \mathsf{pair}(\mathsf{proj}_2(x), \mathsf{proj}_1(x)) : (A \land B) \to (B \land A)$
  \DisplayProof
\]
It does two things at the same time: It shows that $(A \land B) \to (B
\land A)$ has a construction under the BHK interpretation, and that
the lambda term \[\lambda x. \mathsf{pair}(\mathsf{proj}_2 (x),
\mathsf{proj}_1 (x))\] is well typed. In this case, it is a function
which takes arguments that are pairs of type $A \land B$, and returns
pairs of type $B \land A$.

The correspondence goes further. Programs, i.e., lambda terms $t$, can
be executed: we can determine the value they return when we
apply them to an argument~$s$. This is done by reducing the term
$(ts)$ to a normal form, e.g., a term of the form
$\mathsf{proj}_1(\mathsf{pair}(t_1,t_2))$ reduces to~$t_1$. Evaluation
of lambda terms corresponds to normalization of the corresponding
natural deduction proofs (their type derivations). E.g., the reduction
of the term $\mathsf{proj}_1(\mathsf{pair}(t_1,t_2))$ to $t_1$
corresponds to the reduction conversion for $\land$I followed by
$\land$E in natural deduction:
\[
  \AxiomC{$\vdots$}
  \noLine\UnaryInf$\Gamma \fCenter t_1 : A_1$
  \AxiomC{$\vdots$}
\noLine\UnaryInf$\Gamma \fCenter t_2 : A_2$
\RightLabel{$\land$I}
\BinaryInf$\Gamma \fCenter \mathsf{pair}(t_1, t_2) : A_1 \land A_2$
\RightLabel{$\land$E$_1$}
\UnaryInf$\Gamma \fCenter \mathsf{proj}_1 : A_1$
\DisplayProof
\quad\to\quad
\AxiomC{$\vdots$}
\noLine
\UnaryInf$\Gamma \fCenter t_1 : A_1$
\DisplayProof
\]
and the reduction of lambda abstracts applied to terms,
$\beta$-reduction,
\[((\lambda x.t)s) \to t[s/x]\]
corresponds to the reduction conversion for $\to$I followed
by~$\to$E:
\[
  \Axiom$x:A, \Gamma \fCenter x: A$
  \noLine
  \UnaryInfC{$\vdots$}
  \noLine
  \UnaryInf$x:A, \Gamma \fCenter t:B$ 
  \RightLabel{$\to$I}
  \UnaryInf$\Gamma \fCenter \lambda x.t : A \to B$ 
  \AxiomC{$\vdots$}
  \noLine\UnaryInf$\Gamma \fCenter s: A$ 
  \RightLabel{$\to$E}
  \BinaryInf$\Gamma \fCenter ((\lambda x. t)s)B$
  \DisplayProof
\quad\to\quad 
  \AxiomC{$\vdots$} 
  \noLine
  \UnaryInf$\Gamma \fCenter s : A$ 
  \noLine
  \UnaryInfC{$\vdots$}
  \noLine\UnaryInf$\Gamma \fCenter t[s/x] : B$ 
  \DisplayProof 
\]
This correspondence between types and formulas, programs and proofs,
and evaluation and normalization underlies the foundation of typed
programming languages as well as of proof assistants making use of
intuitionistic type theory. 

The lambda calculus we've considered so far as a toy example is not a
useful programming language. To write programs that compute
interesting things, we first have to add data types. This is done by
adding a type $\mathbb{N}$ as a basic type, and terms for the natural
numbers, e.g., $0$, $s(0)$, $s(s(0))$, \dots, together with typing
rules, e.g.,
\[
  \AxiomC{}
  \UnaryInf$\Gamma \fCenter 0:\mathbb{N}$
  \DisplayProof
  \quad
  \Axiom$\Gamma \fCenter t : \mathbb{N}$
  \UnaryInf$\Gamma \fCenter s(t) : \mathbb{N}$
  \DisplayProof
\] 
Programming languages add other basic types, such as boolean values.
We also want to have a way of defining functions by various forms of
recursion.

The power of type systems really only becomes apparent when we add
dependent types and polymorphic functions. In the simply typed lambda
calculus, terms depend on terms (e.g., $\mathsf{proj}_1(x)$ depends
on~$x$) and types depend on types (e.g., $A \land A$ depends on $A$).
A polymorphic term is one that depends not just on terms but also on
types, i.e., one that includes a lambda abstract for types~$\Lambda
X.t$ (of type $\forall X$). This allows us to write terms which work
uniformly on all types. For instance, the term $\Lambda X.\lambda z{:}
X.z$ of type $\forall X(X \to X)$ is a polymorphic identity function:
applied to a type $A$, it returns the identity function of type $A \to
A$.\footnote{It then becomes necessary to mark the type of bound
object variables in the syntax; hence we write $\lambda z{:} X$ and
not just $\lambda z$.} An example of a dependent type might be the
type of $k$-element vectors of natural numbers, $\mathbb{N}^k$. Here
we allow types to depend on terms: The type $\mathbb{N}^t$ depends on
the value of $t$. These features require extensions of the system of
types and of the corresponding lambda calculus, with corresponding
typing rules (and, in the case of dependent types, allowing term
reduction inside of type expressions). 

Logicians have laid the ground work of type systems in the 1970s and
80s. This has had substantial payoffs in the design and implementation
of programming languages. The development of type theories provides
the theoretical foundation of many of the important features of modern
typed programming languages (especially functional languages such as
ML, OCaml, and Haskell). The development of proof systems has also
facilitated the implementation of algorithms to automatically answer
certain questions about programs. E.g., an important aspect of program
development is type checking, i.e., checking whether a program has the
type the programmer claims it to have (in explicitly typed languages,
where the programmer has to provide the type of every function in the
code), or to infer the type a program has (in implicitly typed
languages). This involves finding proofs in the corresponding type
inference calculus. Unification (especially higher-order unification)
plays an important role here too, especially when type inference is
concerned. All of this is true for proof assistants as well. In fact,
the proof assistants based on versions of intuitionistic type theory
discussed in Section~\ref{formal-math} all have their own powerful
programming languages, and proofs produced in them have associated
proof terms, which are essentially programs in these
languages.\footnote{The standard introduction to typed programming
languages is \citet{Pierce2002}. See also \citet{Coquand2022} and
\citet{Zach2019} for discussions aimed at philosophers.}

\section{Conclusion}

We have only scratched the surface of the deep and numerous
connections between logic and mathematics and computer science. The
examples given and episodes described above at least show that the
influence of (philosophical) logic in these fields has been
substantial. Modern foundational mathematical work depends essentially
on logic. Logical methods have yielded results in mainstream
mathematics, especially in algebra. Formalization and verification of
real-life mathematical theorems is no longer just an in-principle
possibility but something actually done on an increasing scale. Large
areas of computer science like automated reasoning using general
purpose and domain-based formalisms, automated verification, type
inference for programming languages, all make essential use of logical
formalisms and methods, which were often pioneered by philosophers. 

The same is true of the reverse: a lot of work in logic done by
mathematicians and computer scientists is also relevant to
philosophical logic. It is also clearly relevant to the philosophies
of mathematics and computer science. Philosophical work on abstraction
principles and theories of truth does not and cannot ignore the
relevant results from mathematical disciplines such as set theory and
proof theory. A lot of recent work in deontic and epistemic logic, in
non-monotonic reasoning, and in belief revision has been done by
computer scientists. Logics motivated by considerations from computer
science, such as Girard's linear logic \citep{Girard1987a} and Parigot's
$\lambda\mu$-calculus \citep{Parigot1992a}, are being studied also by
philosophers working on substructural logics and proof-theoretic
semantics.

Additional connections are just starting to be explored, like those
between type theory as developed by mathematicians and computer
scientists and the use of higher order logic in metaphysics. But it
remains worrisome that the traditional focus in philosophical logic on
what \citet{MartinHjortland2022} call its ``traditional properties''
(e.g., generality, formality, and a prioriticity) excludes large parts
of what goes under the title ``logic'' in mathematics and computer
science. There are of course prominent philosophers of logic in the
anti-exceptionalist camp who reject such a restriction of the field.
\citet{Martin2022,Martin2024} calls for a practice-based approach in
the philosophy of logic, in which anything logicians do counts as
``logic,'' and thus is included in the domain that the philosophy of
logic is concerned with. 

One does not have to adopt such a position wholesale to accept that
some of what goes under the title ``logic'' in mathematics and
computer science is justifiably within the purview of the philosophy
of logic. For instance, higher-order logic and type theories have
sometimes been excluded from the domain of logic proper since the
objects they seem to quantify over are not purely logical, but sets
and set-theoretic functions (most famously by Quine). \emph{Logical}
validity and inference have often been restricted to validity and
inference in (formal regimentations of) natural language. But it is
not just a historical accident that many of our logical systems
(including higher-order logic) were motivated by concerns about
validity and inference in \emph{mathematics}. And as recent
developments in the foundations of mathematics (some discussed above)
have made clear, higher types, functions, and functionals do not
necessarily have to be taken as set-theoretically constructed, but can
be considered primitive notions. In fact, they are so considered by
mathematicians working in foundational systems alternative to set
theory.\footnote{One doesn't have to go to category theory and
intuitionistic type theory to hold this view; even von Neumann, in his
first contributions to the development of set theory considered
functions as more fundamental than sets. And our very notion of a
predicate was introduced as a sui generis function (from objects to
truth values) by Frege.} Perhaps very soon we will look on someone
insisting that type theory is not really logic because it quantifies
over functions like we now do on someone criticizing first-order logic
because it is not just concerned with categorical propositions and
Aristotelian syllogisms.

\paragraph{Acknowledgments.} Thanks to Jeremy Avigad, Steve Awodey,
Matthias Baaz, Chris Fermüller, David Schrittesser and two referees
for Oxford University Press for helpful comments on a previous version.

\bibliography{logic-math-cs}

\end{document}